\documentclass{amsart}

\usepackage{a4,xy,amssymb,diagrams}
\usepackage[hyperindex,pdfmark]{hyperref}

\xyoption{poly}
\xyoption{arc}
\xyoption{knot}
\xyoption{curve}
\xyoption{arrow}
\xyoption{all}

\makeindex

\newcommand{\een}{{\xy 0;/r.08pc/:
(3,-3); (5,-3) **@{-}; (5,5) **@{-}; (0,2) **@{-};
(3,-3); (3,4) **@{-}
\endxy}}

\newcommand{\dummy[1]}{{#1}}
\newcommand{\C}{\mathbb{C}}

\newcommand{\Q}{\mathbb{Q}}
\newcommand{\R}{\mathbb{R}}
\newcommand{\N}{\mathbb{N}}
\newcommand{\Z}{\mathbb{Z}}

\newcommand{\Oscr}{{\mathcal O}}

\newcommand{\wis}[1]{{\text{\em \usefont{OT1}{pag}{m}{n} #1}}}

\providecommand{\abs}[1]{\lvert#1\rvert}

\newtheorem{theorem}{Theorem}[section]
\newtheorem{lemma}[theorem]{Lemma}
\newtheorem{proposition}[theorem]{Proposition}

\theoremstyle{definition}
\newtheorem{definition}[theorem]{Definition}
\newtheorem{example}[theorem]{Example}
\newtheorem{question}[theorem]{Question}

\theoremstyle{remark}

\begin{document}

\sloppy

\title{Necklace Lie algebras and noncommutative Symplectic Geometry}

\author{Raf Bocklandt}
\address{Universiteit Antwerpen (UIA) \\ B-2610 Antwerp (Belgium)}
\email{rbockl@wins.uia.ac.be}
\urladdr{http://win-www.uia.ac.be/u/rbockl/}
\author{Lieven Le Bruyn}
\address{Universiteit Antwerpen (UIA) \\ B-2610 Antwerp (Belgium)}
\email{lebruyn@wins.uia.ac.be}
\urladdr{http://win-www.uia.ac.be/u/lebruyn/}
\thanks{Research Director of the FWO (Belgium)}

\begin{abstract}
Recently, V. Ginzburg proved that Calogero phase space is a coadjoint orbit for
some infinite dimensional Lie algebra coming from noncommutative symplectic
geometry,
\cite{Ginzburg:1999}. In this note we generalize his argument to specific
quotient
varieties of representations of (deformed) preprojective algebras.
This result was also obtained independently by V. Ginzburg \cite{Ginzburg:2000}.
Using results of W. Crawley-Boevey and M. Holland
\cite{CrawleyHolland:1998},
\cite{Crawley:1999} and \cite{Crawley:1999b} we give a combinatorial description
of all the relevant couples $(\alpha,\lambda)$ which are coadjoint orbits.
Finally we explain the coadjoint settings as those for which there is a
Cayley-smooth algebra associated to them.
\end{abstract}

\maketitle

\section{Introduction.}

Over the last couple of years some surprising new results were obtained about the space $Weyl$
of isomorphism classes of right ideals in the first Weyl algebra $A_1(\C)$ linking its structure
to the adelic Grassmannian $Gr^{ad}$ and certain moduli spaces.

Let $\lambda \in \C$, a subset $V \subset \C[x]$ is said to be {\em $\lambda$-primary}
\index{primary!$\lambda$-} if there is some power $r \in \N_+$ such that
\[
(x-\lambda)^r \C[x] \subset V \subset \C[x] \]
A subset $V \subset \C[x]$ is said to be {\em primary decomposable} if it is the finite intersection
\[
V = V_{\lambda_1} \cap \hdots \cap V_{\lambda_r} \]
with $\lambda_i \not= \lambda_j$ if $i \not= j$ and $V_{\lambda_i}$ is a $\lambda_i$-primary
subset. Let $k_{\lambda_i}$ be the codimension of $V_{\lambda_i}$ in $\C[x]$ and consider the
polynomial 
\[
p_V(x) = \prod_{i=1}^r (x-\lambda_i)^{k_{\lambda_i}} \]
Finally, take $W = p_V(x)^{-1} V$, then $W$ is a vectorsubspace of the rational functionfield $\C(x)$
in one variable.

\begin{definition} The {\em adelic Grassmannian} $Gr^{ad}$ is the set of subspaces $W \subset \C(x)$
that arise in this way.
\end{definition}

We can decompose $Gr^{ad}$ in affine cells as follows. For a fixed $\lambda \in \C$ we define
\[
Gr_{\lambda} = \{ W \in Gr^{ad}~\mid~\exists k,l \in \N~:~(x-\lambda)^k \C[x] \subset W \subset
(x-\lambda)^{-l} \C[x] \}
\]
Then, we can write every element $w \in W$ as a {\em Laurent series} \index{Laurent series}
\[
w = \alpha_s (x-\lambda)^s +~\text{higher terms} \]
Consider the increasing set of integers $S = \{ s_0 < s_1 < \hdots \}$ consisting of all {\em degrees}
$s$ of elements $w \in W$. Now, define natural numbers
\[
v_i = i - s_i \qquad \text{then} \qquad v_0 \geq v_1 \geq \hdots \geq v_z = 0 = v_{z+1} = \hdots \]
That is, to $W \in Gr_{\lambda}$ we can associate a {\em partition}
\[
p(W) = (v_0,v_1,\hdots,v_{z-1}) \]
Conversely, if $p$ is a partition of some $n$, then the set of all $W \in Gr_{\lambda}$
with associated partition $p_W = p$ form an affine space $\mathbb{A}^n$ of dimension $n$. Hence,
$Gr_{\lambda}$ has a cellular structure indexed by the set of all partitions.

As $Gr^{ad} = \prod_{\lambda \in \C}' Gr_{\lambda}$ because for any $W \in Gr^{ad}$ there are
uniquely determined $W(\lambda_i) \in Gr_{\lambda_i}$ such that
$W = W(\lambda_1) \cap \hdots \cap W(\lambda_r)$, there is a natural number $n$
associated to $W$ where $n = \abs{p_i}$ where $p_i = p(W(\lambda_i))$ is the partition 
determined by $W(\lambda_i)$. Again, all $W \in Gr^{ad}$ with corresponding
$(\lambda_1,p_1;\hdots;\lambda_r,p_r)$ for an affine cell $\mathbb{A}^n$ of dimension $n$. In his
way, the adelic Grassmannian $Gr^{ad}$ becomes an infinite cellular space with the cells indexed
by $r$-tuples of complex numbers and partitions for all $r \geq 0$.

A surprising connection between
$Gr^{ad}$ and the Calogero system was discovered by G. Wilson \index{Wilson, G.} in \cite{Wilson:1998}.

\begin{theorem} Let $Gr^{ad}(n)$ be the collection of all cells of dimension $n$ in $Gr^{ad}$, then
there is a set-theoretic bijection
\[
Gr^{ad}(n) \longleftrightarrow Calo_n \]
between $Gr^{ad}(n)$ and the phase space of $n$ Calogero particles, that is, the orbit space of
couples of $n \times n$ matrices $(X,Y)$ such that the rank of $[X,Y] - \een_n$ is one under
simultaneous conjugation.
\end{theorem}

The connection between right ideals of $A_1(\C)$ and $gr^{ad}$ is contained (in disguise) in
the paper of R. Cannings \index{Cannings, R.} and M. Holland \index{Holland, M.} \cite{CanningsHolland:1994}.
$A_1(\C)$ acts as differential operators on $\C[x]$ and for every right ideal $I$ of $A_1(\C)$ they show
that $I.\C[x]$ is primary decomposable. Conversely, if $V \subset \C[x]$ is primary decomposable, they
associate the right ideal
\[
I_V = \{ \theta \in A_1(\C)~\mid~\theta.\C[x] \subset V \}
\]
of $A_1(\C)$ to it. Moreover, isomorphism classes of right ideals correspond to studying
primary decomposable subspaces under multiplication with polynomials. Hence,
\[
Gr^{ad} \simeq Weyl \]
The group $Aut~A_1(\C)$ of $\C$-algebra automorphisms of $A_1(\C)$ acts on the set of right
ideals of $A_1(\C)$ and respects the notion of isomorphism whence acts on $Weyl$. 
The group $Aut~A_1(\C)$ is generated by automorphisms $\sigma_i^f$ defined by
\[
\begin{cases}
\sigma_1^f(x) &= x + f(y) \\
\sigma_1^f(y) &= y \end{cases} \quad \text{with $f \in \C[y]$,} \quad 
\begin{cases}
\sigma_2^f(x) &= x  \\
\sigma_2^f(y) &= y +f(x) \end{cases} \quad \text{with $f \in \C[x]$} \]
For a natural number $n \geq 1$ we define the {\em $n$-th canonical right ideal} 
\index{canonical right ideal} of $A_1(\C)$ to be
\[
\mathfrak{p}_n = x^{n+1} A_1(\C) + (xy+n) A_1(\C) . \]
One can show that $\mathfrak{p}_n \not\simeq \mathfrak{p}_m$ whenever $n \not= m$ so the
isomorphism classes $[ \mathfrak{p}_n ]$ are distinct points in $Weyl$ for all $n$. We 
define
\[
Weyl_n = Aut~A_1(\C) . [ \mathfrak{p}_n ] = \{~[ \sigma(\mathfrak{p}_n) ] \quad
\forall \sigma \in Aut~A_1(\C) \} \]
the orbit in $Weyl$ of the point $[ \mathfrak{p}_n ]$ under the action of the automorphism
group.

Yu. Berest \index{Berest, Y.} and G. Wilson \index{Wilson, G.} proved in \cite{BerestWilson:1999}
that the Cannings-Holland correspondence respects the automorphism orbit decomposition.

\begin{theorem} We have $Weyl = \bigsqcup_n~Weyl_n$ and there are set-theoretic bijections
\[
Weyl_n \longleftrightarrow Gr^{ad}(n) \]
whence also with $Calo_n$.
\end{theorem}

An early indication that a decomposition into moduli spaces might be possible can be found in
the paper \cite{LeBruyn:1995} of the second author. Recently, a similar idea was pursued by
A. Kapustin, A. Kutzetsov and D. Orlov \cite{Kapustin}.

If we trace the action of $Aut~A_1(\C)$ on $Weyl_n$ through all the identifications, we get a
transitive action of $Aut~A_1(\C)$ on $Calo_n$. However, this action is non-differentiable
hence highly non-algebraic. Berest and Wilson asked whether it is possible to identify $Calo_n$
with a coadjoint orbit in some infinite dimensional Lie algebra.

This conjecture was proved by V. Ginzburg \cite{Ginzburg:1999} using noncommutative symplectic
geometry as sketched by M. Kontsevich \cite{Kontsevich:1993}. After reading his preprint it
became clear to us that his method could be used almost verbatim for the quotient varieties
of representations of deformed preprojective algebras. The two crucial steps are  (1) 
invariants of quiver representations are generated by traces along oriented cycles, proved in
\cite{LeBruynProcesi:1990} and (2) acyclicity of the relative (!) noncommutative deRham cohomology
for path algebras of quivers.

In the first sections of this paper we carry out the second project in some detail. When applied
to the double $\Q$ of a quiver, the noncommutative functions acquire the structure of an infinite
dimensional Lie algebra $\N_Q$ which we call the necklace Lie algebra of a quiver. A major result
asserts that this Lie algebra is a central extension of the Lie algebra of symplectic derivations,
that is the Lie algebra corresponding to the vertex-preserving automorphisms of the path algebra
$\C \Q$ preserving the moment element $m = \sum_a [a,a^*]$.

Recall that the deformed preprojective algebra $\Pi_{\lambda}$ is defined to be the quotient of
$\C \Q$ by the twosided ideal generated by $m - \lambda$. Our generalization of Ginzburg's result
on the coadjointness of Calogero space can then be stated as :

\begin{theorem} The variety $iss_{\alpha}~\Pi_{\lambda}$ of isomorphism classes of semisimple
$\alpha$-dimensional representations of the deformed preprojective algebra $\Pi_{\lambda}$ is
a coadjoint orbit of the necklace Lie algebra $\N_Q$ whenever $\alpha$ is a minimal non-zero
element of $\Sigma_{\lambda}$ the set of dimension vectors of simple representations of $\Pi_{\lambda}$.
\end{theorem}

W. Crawley-Boevey has given a combinatorial description of the set $\Sigma_{\lambda}$ in \cite{Crawley:1999}.
We recover the Calogero case back and prove that the spaces appearing in a conjectural extension to
arbitrary extended Dynkin quivers (conjectured by M. Holland and W. Crawley-Boevey) are all
coadjoint orbits.

In the last two sections we try to explain why precisely these $(\alpha,\lambda)$ couples appear from
the viewpoint of noncommutative geometry. As the path algebra $\C \Q$ is a formally smooth algebra
as in \cite{CuntzQuillen:1995}, its representation spaces $rep_{\alpha}~\Q$ are smooth varieties and
the noncommutative functions and differential forms induce invariant classical functions and forms
on these varieties and their quotient varieties. On the other hand, we will show that the deformed
preprojective algebra $\Pi_{\lambda}$ is {\em not} formally smooth and so should be viewed as a
singular subvariety of the noncommutative manifold corresponding to $\C \Q$. As such, the
noncommutative vectorfields on $\C \Q$ (the Lie algebra $\N_Q$) have rather unpredictable
behavior on the singular closed subvariety $\Pi_{\lambda}$. However, for those dimension vectors
$\alpha$ such that $rep_{\alpha}~\Pi_{\lambda}$ is smooth (that is, $\Pi_{\lambda}$ is an
$\alpha$-smooth subvariety of $\C \Q$) things should work out. We conjecture that the corresponding
$\alpha$ are precisely the minimal elements in $\Sigma_{\lambda}$ (the coadjoint orbits). We prove this
for the preprojective algebra $\Pi_0$ and prove a variation (using hyper-K\"ahler reductions) for
$\Pi_{\lambda}$. These results are based on the calculation of $Ext^1$'s of representations of
$\Pi_0$ due to W. Crawley-Boevey \cite{Crawley:1999b}.

\section{Necklace Lie algebras.}

In this section we introduce the main object of this note in a purely
combinatorial
way. Recall that a {\em quiver} $Q$ is a finite directed graph on a set of
vertices
$Q_v = \{ v_1,\hdots,v_k \}$, having a finite set $Q_a = \{ a_1,\hdots,a_l
\}$ of arrows,
where we allow loops as well as multiple arrows between vertices. An arrow $a$
with starting vertex
$s(a) = v_i$ and terminating vertex $t(a) = v_j$ will be depicted as
$\xy
\POS (15,0) *\cir<4pt>{}*+{\txt\tiny{i}} ="v2"
   , (0,0) *\cir<4pt>{}*+{\txt\tiny{j}} ="v3"
\POS"v2" \ar "v3"_{a} \endxy$. The quiver information is encoded in the
{\em Euler form}
which is the bilinear form on $\Z^k$ determined by the matrix $\chi_Q \in
M_k(\Z)$ with
\[
\chi_{ij} = \delta_{ij} - \#~\{~a \in Q_a~\mid~\xy
\POS (15,0) *\cir<4pt>{}*+{\txt\tiny{i}} ="v2"
   , (0,0) *\cir<4pt>{}*+{\txt\tiny{j}} ="v3"
\POS"v2" \ar "v3"_{a} \endxy~\}
\]
The symmetrization $T_Q = \chi_Q + \chi_Q^{tr}$ of this matrix determines
the {\em Tits form} of the
quiver $Q$. An oriented cycle $c = a_{i_u} \hdots a_{i_1}$ of length $u
\geq 1$ is a concatenation of
arrows in $Q$ such that $t(a_{i_j}) = s(a_{i_{j+1}})$ and $t(a_{i_u}) =
s(a_{i_1})$. In
addition to these there are $k$ oriented cycles
$e_i$ of length $0$ corresponding to the vertices of $Q$. All oriented
cycles $c'$ obtained from $c$
by cyclically permuting the arrow components are said to be equivalent to
$c$. A {\em necklace word}
$w$ for $Q$ is an equivalence class of oriented cycles in the quiver $Q$.

The {\em double quiver} $\Q$ of $Q$ is the quiver obtained by adjoining to
every arrow (or loop)
$\xy
\POS (15,0) *\cir<4pt>{}*+{\txt\tiny{i}} ="v2"
   , (0,0) *\cir<4pt>{}*+{\txt\tiny{j}} ="v3"
\POS"v2" \ar "v3"_{a} \endxy$ in $Q$ an arrow in the opposite direction
$\xy
\POS (15,0) *\cir<4pt>{}*+{\txt\tiny{i}} ="v2"
   , (0,0) *\cir<4pt>{}*+{\txt\tiny{j}} ="v3"
\POS"v3" \ar "v2"^{a^*} \endxy$. That is,
$
\chi_{\Q} = T_Q - \een_k$.

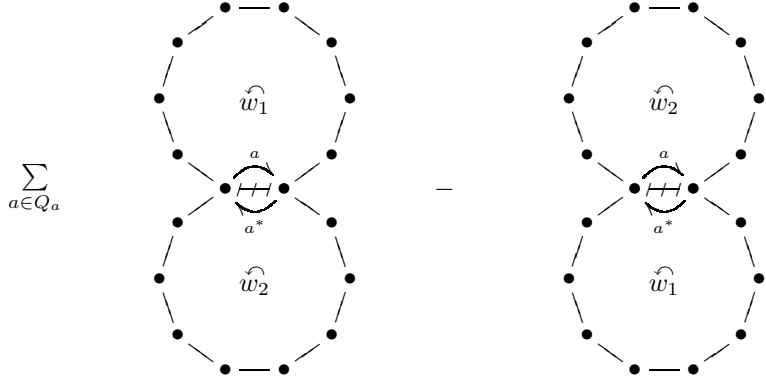
\begin{figure}
\[
\begin{xy}/r3pc/:
{\xypolygon10{~*{\bullet}~>>{}}},
"1" *+{\txt{\tiny $u$}},"8"="a","10"="a1","9"="a2","0"="c1",
"c1"+(0,-1.9),
{\xypolygon10{~*{\bullet}~>>{}}},
"5" *+{\txt{\tiny $v$}},"4"="b","6"="b1","5"="b2","0"="c2",
"a";"a2" **@{/};
"a" +(-2,0) *+{\txt{$\underset{a \in Q_a}{\sum}$}};
"c1" *+{\txt{$\overset{\curvearrowleft}{w_1}$}};
"c2" *+{\txt{$\overset{\curvearrowleft}{w_2}$}};
\POS"a" \ar@/^2ex/^{\txt{\tiny{$a$}}} "a2"
\POS"a2" \ar@/^2ex/^{\txt{\tiny{$a^*$}}} "a"
\end{xy}~\qquad
 \begin{xy}/r3pc/:
{\xypolygon10{~*{\bullet}~>>{}}},
"1" *+{\txt{\tiny $u$}},"8"="a","10"="a1","9"="a2","0"="c1",
"c1"+(0,-1.9),
{\xypolygon10{~*{\bullet}~>>{}}},
"5" *+{\txt{\tiny $v$}},"4"="b","6"="b1","5"="b2","0"="c2",
"a";"a2" **@{/};
"a"+(-2,0) *+{-};
"c1" *+{\txt{$\overset{\curvearrowleft}{w_2}$}};
"c2" *+{\txt{$\overset{\curvearrowleft}{w_1}$}};
\POS"a" \ar@/^2ex/^{\txt{\tiny{$a$}}} "a2"
\POS"a2" \ar@/^2ex/^{\txt{\tiny{$a^*$}}} "a"
\end{xy}
\]
\caption{Lie bracket $[ w_1,w_2 ]$ in $\N_Q$.}
\label{bracket}
\end{figure}

The {\em necklace Lie algebra} $\N_Q$ for the quiver $Q$ has as  basis the set
of all necklace words $w$ for the {\em double} quiver $\Q$ and where the
Lie bracket $[w_1,w_2]$
is determined as in figure~\ref{bracket}.
That is, for every arrow $a \in Q_a$ we look for an occurrence of $a$ in
$w_1$ and of $a^*$ in
$w_2$. We then open up the necklaces by removing these factors and regluing
the open ends together
to form a new necklace word. We repeat this operation for {\em all}
occurrences of $a$ (in $w_1$)
and $a^*$ (in $w_2$). We then replace the roles of $a^*$ and $a$ and redo
this operation with a
minus sign. Finally, we add up all these obtained necklace words for all
arrows $a \in Q_a$. Using
this graphical description it is a pleasant exercise to verify the Jacobi
identity for $\N_Q$.

\section{An acyclicity result.}

The {\em path algebra} $\C Q$ of a quiver $Q$ has as basis the set of all
oriented paths
$p = a_{i_u} \hdots a_{i_1}$ of length $u \geq 1$ in the quiver, that is
$s(a_{i_{j+1}}) = t(a_{i_j})$
together with the vertex-idempotents $e_i$ of length zero. Multiplication in
$\C Q$ is induced by (left) concatenation of paths. More precisely, $1= e_1
+ \hdots + e_k$ is a
decomposition of $1$ into mutually orthogonal idempotents and further we define
\begin{itemize}
\item{$e_j.a$ is always zero unless \ \  $\xy
\POS (0,0) *\cir<4pt>{}*+{\txt\tiny{j}} ="v2"
   , (15,0) *\cir<4pt>{} ="v3"
\POS"v3" \ar "v2"_{a} \endxy$ \ \ in which case it is the path $a$,}
\item{$a.e_i$ is always zero unless \ \ $\xy
\POS (15,0) *\cir<4pt>{}*+{\txt\tiny{i}} ="v2"
   , (0,0) *\cir<4pt>{} ="v3"
\POS"v2" \ar "v3"_{a} \endxy$ \ \ in which case it is the path $a$,}
\item{$a_i.a_j$ is always zero unless \ \ $\xy
\POS (15,0) *\cir<4pt>{} ="v2"
   , (0,0) *\cir<4pt>{} ="v1"
   , (30,0) *\cir<4pt>{} ="v3"
\POS"v2" \ar "v1"_{a_i}
\POS"v3" \ar "v2"_{a_j} \endxy$ \ \ in which case it is the path $a_ia_j$.}
\end{itemize}
Path algebras of quivers are the archetypical examples of {\em formally
smooth algebras} as
introduced and studied in \cite{CuntzQuillen:1995}.

In this section we will generalize Kontsevich's acyclicity result for the
noncommutative
deRham cohomology of the free algebra \cite{Kontsevich:1993} to that of the
path algebra $\C Q$.
The crucial idea is to consider the {\em relative} differential forms
(as defined in \cite{CuntzQuillen:1995}) of $\C Q$ with respect to the
semisimple subalgebra $V = \C \times \hdots \times \C$ generated by the
vertex idempotents. The idea being that in considering quiver representations
one works in the category of $V$-algebras rather than $\C$-algebras.

For a subalgebra $B$ of $A$, let $\overline{A}_B$ denote the cokernel of the
inclusion as $B$-bimodule.
The space of relative differential forms of degree $n$ of $A$ with respect
to $B$ is
\[
\Omega^n_B~A = A \otimes_B \underbrace{\overline{A}_B \otimes_B \hdots
\otimes_B \overline{A}_B}_n \]
The space $\Omega^{\bullet}_B~A$ is given a differential graded algebra
structure by taking the multiplication
\[
(a_0,\hdots,a_n)(a_{n+1},\hdots,a_m) = \sum_{i=0}^n
(-1)^{n-i}(a_0,\hdots,a_{i-1},a_ia_{i+1},a_{i+2},\hdots,a_m)
\]
and the differential $d(a_0,\hdots,a_n) = (1,a_0,\hdots,a_n)$, see
\cite{CuntzQuillen:1995}. Here,
$(a_0,\hdots,a_n)$ is a representant of the class $a_0 da_1 \hdots da_n \in
\Omega^n_B~A$ and we recall that
$\Omega^{\bullet}_B~A$ s generated by the $a$ and $da$ for all $a \in A$.
The {\em relative
cohomology} $H^n_B~A$ is defined as the cohomology of the complex
$\Omega^{\bullet}_B~A$.

For $\theta \in Der_{B}~A$, the Lie algebra
of $B$-derivations of $A$ (that is $\theta$ is a derivation of $A$ and
$\theta(B) = 0$), we define a
degree preserving derivation $L_{\theta}$ and a degree $-1$
super-derivation $i_{\theta}$ on
$\Omega^{\bullet}_B~A$
\[
\xy
\POS (10,0) *+{\txt{$\Omega^{n-1}_B~A$}} ="a" ,
\POS (40,0) *+{\txt{$\Omega^{n}_B~A$}} ="b" ,
\POS (70,0) *+{\txt{$\Omega^{n+1}_B~A$}} ="c"
 \POS"a" \ar@(dl,dr)_{\txt{$L_{\theta}$}}
 \POS"b" \ar@(dl,dr)_{\txt{$L_{\theta}$}}
 \POS"c" \ar@(dl,dr)_{\txt{$L_{\theta}$}}
\POS"a" \ar @/^4ex/^{\txt{$d$}} "b"
\POS"b" \ar @/^4ex/^{\txt{$i_{\theta}$}} "a"
\POS"b" \ar @/^4ex/^{\txt{$d$}} "c"
\POS"c" \ar @/^4ex/^{\txt{$i_{\theta}$}} "b"
\endxy
\]
by the rules
\[
\begin{cases}
L_{\theta}(a) = \theta(a) & \quad L_{\theta}(da) = d~\theta(a) \\
i_{\theta}(a) = 0 & \quad i_{\theta}(da) = \theta(a)
\end{cases}
\]
for all $a \in A$. We have the Cartan homotopy formula
$
L_{\theta} = i_{\theta} \circ d + d \circ i_{\theta}$
as both sides are degree preserving derivations on $\Omega^{\bullet}_B~A$
and they agree on all the
generators $a$ and $da$ for $a \in A$.

\begin{lemma} \label{contraction} Let $\theta, \gamma \in Der_{B}~A$, then
we have on $\Omega^{\bullet}_B~A$
the identities of operators
\[
\begin{cases}
L_{\theta} \circ i_{\gamma} - i_{\gamma} \circ L_{\theta} = [L_{\theta},
i_{\gamma} ] &= i_{[\theta,
\gamma]} = i_{\theta \circ \gamma - \gamma \circ \theta} \\
L_{\theta} \circ L_{\gamma} - L_{\gamma} \circ L_{\theta} = [L_{\theta},
L_{\gamma} ] &=
L_{[\theta,\gamma]} = L_{\theta \circ \gamma - \gamma \circ \theta}
\end{cases}
\]
\end{lemma}

\begin{proof} Consider the first identity. By definition both sides are
degree $-1$ super-derivations
on $\Omega^{\bullet}_B~A$ so it suffices to check that they agree on
generators. Clearly, both sides give $0$
when evaluated on $a \in A$ and for $da$ we have
\[
(L_{\theta} \circ i_{\gamma} - i_{\gamma} \circ L_{\theta}) da =
L_{\theta}~\gamma(a) - i_{\gamma}~d~\theta(a) = \theta~\gamma(a) -
\gamma~\theta(a) = i_{[\theta,\gamma]}(da)
\]
A similar argument proves the second identity.
\end{proof}

Specialize to the quiver-case with $A = \C Q$ the path algebra and $B = V =
\C^k$ the vertex algebra.

\begin{lemma} Let $Q$ be a quiver on $k$ vertices, then a basis for
$\Omega^n_V~\C Q$ is
given by the elements
\[
p_0 d p_1 \hdots d p_n \]
where $p_i$ is an oriented path in the quiver such that $length~p_0 \geq 0$ and
$length~p_i \geq 1$ for $1 \leq i \leq n$ and such that the starting point
of $p_i$ is the endpoint
of $p_{i+1}$ for all $1 \leq i \leq n-1$.
\end{lemma}

\begin{proof} Clearly $l(p_i) \geq 1$ when $i \geq 1$ or $p_i$ would be a
vertex-idempotent whence
in $V$. Let $v$ be the starting point of $p_i$ and $w$ the end point of
$p_{i+1}$ and assume that
$v \not= w$, then
\[
p_i \otimes_V p_{i+1} = p_i v \otimes_V w p_{i+1} = p_i vw \otimes_V
p_{i+1} = 0
\]
from which the assertion follows.
\end{proof}

\begin{proposition} Let $Q$ be a quiver on $k$ vertices, then the relative
differential form-algebra
 $\Omega_V^{\bullet}~\C Q$
is formal. In fact, the complex is acyclic
\[
\begin{cases}
H^0_V~\C Q &\simeq \C \times \hdots \times \C~\text{($k$ factors)} \\
H^n_V~\C Q &\simeq 0 \qquad \forall n \geq 1
\end{cases}
\]
\end{proposition}

\begin{proof} Define the {\em Euler derivation}
$E$ on $\C Q$ by the rules that
\[
E(e_i) = 0~\forall~1 \leq i \leq k \qquad \text{and} \qquad E(a) =
a~\forall a \in Q_a \]
By induction on the length $l(p)$ of an oriented path $p$ in the quiver $Q$
one easily verifies that
$E(p) = l(p) p$. By induction one can also proof that
$L_E(p_0dp_1\dots dp_n)= l(p_0) + \cdots +l(p_n)$.
This implies that $L_E$ is a bijection on each $\Omega_V^i~\C Q$, where $i>1$ and
on $\Omega_V^0~\C Q$, $L_E$ has V as its kernel.
By applying the Cartan homotopy formula for $L_E$, we obtain that the 
complex is acyclic.
\end{proof}

The complex $\Omega_V^{\bullet}~\C Q$ induces the {\em relative Karoubi
complex}
\[
\wis{dR}^0_V~\C Q \rTo^d \wis{dR}^1_V~\C Q \rTo^d \wis{dR}^2_V~\C Q \rTo^d
\hdots \]
with
\[
\wis{dR}^n_V~\C Q = \dfrac{\Omega^n_V~\C Q}{\sum_{i=0}^n [~\Omega^i_V~\C Q,
\Omega^{n-i}_V~\C Q~]} \]
In this expression the brackets denote supercommutators with respect to
the grading on $\Omega_V^{\bullet}~\C Q$.
In the commutative case, $\wis{dR}^0$ are the functions on the manifold and
$\wis{dR}^1$ the $1$-forms.

\begin{lemma} A $\C$-basis for the noncommutative functions
\[
\wis{dR}^0_V~\C Q \simeq \dfrac{\C Q}{[~\C Q,\C Q~]} \]
are the necklace words in the quiver $Q$.
\end{lemma}

\begin{proof}
Let $\mathbb{W}$ be the $\C$-space spanned by all necklace words $w$ in $Q$
and define a linear map
\[
\C Q \rOnto^n \mathbb{W} \qquad \begin{cases} p \mapsto w_p
&\qquad~\text{if $p$ is a cycle} \\
p \mapsto 0 &\qquad~\text{if $p$ is not} \end{cases} \]
for all oriented paths $p$ in the quiver $Q$, where $w_p$ is the necklace
word in $Q$ determined by
the oriented cycle $p$. Because $w_{p_1p_2} = w_{p_2p_1}$ it follows that
the commutator subspace
$[\C Q,\C Q]$ belongs to the kernel of this map. Conversely, let
\[
x = x_0 + x_1 + \hdots + x_m \]
be in the kernel where $x_0$ is a linear combination of non-cyclic paths
and $x_i$ for $1 \leq i \leq m$
is a linear combination of cyclic paths mapping to the same necklace word
$w_i$, then $n(x_i) = 0$
for all $i \geq 0$. Clearly, $x_0 \in [\C Q,\C Q]$ as we can write every
noncyclic
path $p = a.p' = a.p'-p'.a$ as a commutator. If $x_i = a_1p_1 + a_2p_2 +
\hdots + a_lp_l$ with
$n(p_i) = w_i$, then $p_1 = q.q'$ and $p_2 = q'.q$ for some paths $q,q'$
whence $p_1-p_2$ is a commutator.
But then, $x_i = a_1(p_1-p_2) + (a_2-a_1)p_2 + \hdots + a_l p_l$ is a sum
of a commutator and a linear
combination of strictly fewer elements. By induction, this shows that $x_i
\in [\C Q,\C Q]$.
\end{proof}

\begin{lemma} $\wis{dR}^1_V~\C Q$ is isomorphic as $\C$-space to
\[
\underset{\xy
 \POS (-10,0) *+{\txt{\tiny j}}*\cir<6pt>{} ="a" , (0,0) *+{\txt{\tiny
i}}*\cir<6pt>{}="b"
 \POS"b" \ar@(l,r)^{\txt{\tiny a}} "a"\endxy}{\bigoplus}~v_i.\C Q.v_j~da =
\underset{\xy
 \POS (-10,0) *+{\txt{\tiny j}}*\cir<6pt>{} ="a" , (0,0) *+{\txt{\tiny
i}}*\cir<6pt>{}="b"
 \POS"b" \ar@(l,r)^{\txt{\tiny a}} "a"\endxy}{\bigoplus}~\xy \POS (-15,0)
*+{\txt{\tiny i}}*\cir<6pt>{} ="a" , (0,0) *+{\txt{\tiny j}}*\cir<6pt>{}="b"
 \POS"b" \ar@{.>}@(ul,ur) "a"
 \endxy~d~\xy
 \POS (-10,0) *+{\txt{\tiny j}}*\cir<6pt>{} ="a" , (0,0) *+{\txt{\tiny
i}}*\cir<6pt>{}="b"
 \POS"b" \ar@(l,r)^{\txt{\tiny a}} "a"\endxy
\]
\end{lemma}

\begin{proof}
If $p.q$ is not a cycle, then $p dq = [p,dq]$ and so vanishes in
$\wis{dR}^1_V~\C Q$ so we
only have to consider terms $p dq$ with $p.q$ an oriented cycle in $Q$. For
any three paths
$p,q$ and $r$ in $Q$ we have the equality
\[
[p.q dr] = pq dr - qd(rp) + qr dp \]
whence in $\wis{dR}^1_V~\C Q$ we have relations allowing to reduce the
length of the differential part
\[
q d(rp) = p q dr + q r dp \]
so $\wis{dR}^1_V~\C Q$ is spanned by terms of the form $p da$ with $a \in
Q_a$ and $p.a$ an
oriented cycle in $Q$. Therefore, we have a surjection
\[
\Omega^1_V~\C Q \rOnto \underset{\xy
 \POS (-10,0) *+{\txt{\tiny j}}*\cir<6pt>{} ="a" , (0,0) *+{\txt{\tiny
i}}*\cir<6pt>{}="b"
 \POS"b" \ar@(l,r)^{\txt{\tiny a}} "a"\endxy}{\bigoplus}~v_i.\C Q.v_j~da
 \]
By construction, it is clear that $[\Omega^0_V~\C Q,\Omega^1_{rel}~\C Q]$
lies in the kernel of
this map and using an argument as in the lemma above one shows also the
converse inclusion.
\end{proof}

Using the above descriptions of $\wis{dR}^i_V~\C Q$ for $i=0,1$ and the
differential
$\wis{dR}^0_V~\C Q \rTo^d \wis{dR}^1_V~\C Q$ we can define {\em partial
differential
operators} associated to any arrow $\xy
 \POS (-10,0) *+{\txt{\tiny j}}*\cir<6pt>{} ="a" , (0,0) *+{\txt{\tiny
i}}*\cir<6pt>{}="b"
 \POS"b" \ar@(l,r)^{\txt{\tiny a}} "a"\endxy$ in $Q$.
\[
\frac{\partial}{\partial a}~:~\wis{dR}^0_{rel}~\C Q \rTo v_i \C Q v_j \qquad
\text{by} \qquad df = \sum_{a \in Q_a}~\frac{\partial f}{\partial a} d a \]
To take the partial derivative of a necklace word $w$ with respect to an
arrow $a$, we run through $w$ and each time we encounter $a$ we open the
necklace by removing
that occurrence of $a$ and then take the sum of all the paths obtained.

Defining the {\em relative deRham cohomology} $H^n_{dR}~\C Q$ to be the
cohomology of the Karoubi
complex and observing that the operators $L_{\theta}$ and $i_{\theta}$ on
$\Omega^{\bullet}_V~\C Q$
induce operators on the Karoubi complex, we have the {\em acyclicity result}

\begin{theorem} The relative Karoubi complex is acyclic. In particular,
\[
\begin{cases}
H^0_{dR}~\C Q &\simeq~V \\
H^n_{dR}~\C Q &\simeq~0 \qquad \forall n \geq 1
\end{cases}
\]
\end{theorem}

\section{Symplectic interpretation.}

In this section we use the acyclicity result to give a Poisson
interpretation to the Lie bracket
in $\N_Q$. This generalizes the {\em Kontsevich bracket}
\cite{Kontsevich:1993} in the free case
to path algebras of doubles of quivers. If $Q$ is a quiver with double
quiver $\Q$, then we can
define a canonical {\em symplectic structure} on the path algebra of the
double $\C \Q$ determined
by the element
\[
\omega = \sum_{a \in Q_a} da^* da \in \wis{dR}^2_V~\C \Q \]
As in the commutative case, $\omega$ defines a bijection between the
noncommutative $1$-forms
$\wis{dR}^1_V~\C \Q$ and the {\em noncommutative vectorfields} which are
defined to be the
$V$-derivations of $\C \Q$. This correspondence is
\[
Der_V~\C \Q \rTo^{\tau} \wis{dR}^1_V~\C \Q \qquad \text{given by} \qquad
\tau(\theta) = i_{\theta}(\omega)
\]
In analogy with the commutative case we define a derivation $\theta \in
Der_V~\C \Q$ to be
{\em symplectic} if and only if $L_{\theta} \omega = 0 \in \wis{dR}^2_V~\C \Q$
and denote the subspace of symplectic derivations by
$Der_{\omega}~\C \Q$. It follows from the homotopy formula
and the fact that $\omega$ is a closed form, that $\theta \in
Der_{\omega}~\C \Q$ implies
$
L_{\theta} \omega = d i_{\theta} \omega = d\tau(\theta) = 0$.
That is, $\tau(\theta)$ is a closed form which by the acyclicity of the
Karoubi complex shows that
it must be an exact form. That is we have an isomorphism of exact sequences
of $\C$-vectorspaces
\[
\begin{diagram}
0 & \rTo & V & \rTo & \wis{dR}^0_V~\C \Q & \rTo^d & (\wis{dR}^1_V~\C
\Q)_{exact} & \rTo & 0 \\
& & \dTo^{=} & & \dTo^{\simeq} & & \dTo^{\tau^{-1}} & & \\
0 & \rTo & V & \rTo & \frac{\C \Q}{[\C \Q,\C \Q]} & \rTo & Der_{\omega}~\C
\Q & \rTo & 0
\end{diagram}
\]
The symplectic structure $\omega$ defines a Poisson bracket on the
noncommutative functions.

\begin{definition} Let $Q$ be a quiver and $\Q$ its double. The {\em
Kontsevich bracket} on the necklace words
in $\Q$, $\wis{dR}^0_V~\C \Q$ is defined to be
\[
\{ w_1,w_2 \}_K = \sum_{a \in Q_a}~(\frac{\partial w_1}{\partial a}
\frac{\partial w_2}{\partial a^*}
- \frac{\partial w_1}{\partial a^*} \frac{\partial w_2}{\partial
a})~\wis{mod}~[ \C \Q,\C \Q ]
\]
By the description of the partial differential operators it is clear that
$\wis{dR}^0_V~\C \Q$
with this bracket is isomorphic to the necklace Lie algebra $\N_Q$.
\end{definition}

The symplectic derivations $Der_{\omega}~\C \Q$ have a natural Lie algebra
structure by commutators
of derivations. We will show that $\tau^{-1} \circ d$ is a Lie algebra
morphism.

For every necklace word $w$ we have a symplectic derivation $\theta_w = \tau^{-1} d w$
defined by
\[
\begin{cases}
\theta_w(a) &= -\frac{\partial w}{\partial a^*} \\
\theta_w(a^*) &=  \frac{\partial w}{\partial a}
\end{cases}
\]
With this notation we get the following interpretations of the Kontsevich
bracket
\[
\{ w_1,w_2 \}_K = i_{\theta_{w_1}}(i_{\theta_{w_2}} \omega) =
L_{\theta_{w_1}}(w_2) = - L_{\theta_{w_2}}(w_1) \]
where the next to last equality follows because $i_{\theta_{w_2}} \omega =
d w_2$ and the fact that
$i_{\theta_{w_1}}(d w) = L_{\theta_{w_1}}(w)$ for any $w$. More generally,
for any $V$-derivation
$\theta$ and any necklace word $w$ we have the equation
\[
i_{\theta}(i_{\theta_w} \omega) = L_{\theta}(w). \]

When we look at the image of the Kontsevich bracket under $\tau^{-1}d$, we 
obtain the following
\[
\begin{split}
\tau^{-1}d\{w_1,w_2\}_K &= \tau^{-1}dL_{\theta_{w_1}}w_2\\
&= \tau^{-1}L_{\theta_{w_1}}dw_2\\
&= \tau^{-1}L_{\theta_{w_1}}i_{\theta_{w_2}}\omega\\
&= \tau^{-1}([L_{\theta_{w_1}},i_{\theta_{w_2}}]+i_{\theta_{w_2}}L_{\theta_{w_1}})\omega\\
&=\tau^{-1}i_{[{\theta_{w_1}},{\theta_{w_2}}]}\omega
\\
&=[{\theta_{w_1}},{\theta_{w_2}}]
\end{split}
\]
Above we made use of the fact that $L_\theta$ commutes with $d$, and the defining
equation $dw_2= i_{\theta_{w_2}}\omega$. 
In the fourth line we omitted the last term because $\theta_{w_1}$ is a symplectic derivation.
Finally lemma \ref{contraction} enabled us to transform the commutator in $i$ and $L$ 
to of commutator of the derivations $\theta_{w_1}$ and $\theta_{w_2}$. This calculation concluded the proof of :

\begin{theorem} With notations as before, $\wis{dR}^0_{rel}~\C Q^d$ with the
Kontsevich bracket is isomorphic to the necklace Lie algebra $\N_Q$, and
the sequence
\[
0 \rTo V \rTo \N_Q \rTo^{\tau^{-1} d} Der_{\omega}~\C \Q \rTo 0 \]
is an exact sequence (hence a central extension) of Lie algebras.
\end{theorem}

\section{Coadjoint orbits.}

Consider a dimension vector $\alpha = (n_1,\hdots,n_k)$, that is, a
$k$-tuple of natural numbers, then
the space of $\alpha$-dimensional representations of the double quiver
$\Q$, $rep_{\alpha}~\Q$ can be
identified via the trace pairing with the cotangent bundle
$T^*~rep_{\alpha}~Q$ of the space of
$\alpha$-dimensional representations of the quiver $Q$, see for example
\cite{Crawley:1999}, and as
such acquires a natural symplectic structure. The natural action of the
basechange group
$GL(\alpha) = GL_{n_1} \times \hdots \times GL_{n_k}$ on $rep_{\alpha}~\Q$
is symplectic and induces
a Poisson structure on the coordinate ring as well as on the ring of
polynomial quiver invariants,
which are generated by traces along oriented cycles by
\cite{LeBruynProcesi:1990}.

The symplectic derivations $Der_{\omega}~\C \Q$ correspond to the 
$V$-automorphisms of the path algebra of the double $\C \Q$ preserving the {\em moment element}
\[
m = \sum_{a \in Q_a} [a,a^*] \in \C \Q \]
For this reason it is natural to consider the {\em complex moment map}
\[
rep_{\alpha}~\Q \rTo^{\mu_{\C}} M_{\alpha}^0(\C) \qquad V \mapsto \sum_{a
\in Q_a} [V_a,V_{a^*}]
\]
where $M_{\alpha}^0(\C)$ is the subspace of $k$-tuples $(m_1,\hdots,m_k)
\in M_{n_1}(\C) \oplus \hdots
\oplus M_{n_k}(\C)$ such that $\sum_i tr(m_i) = 0$, that is
$M_{\alpha}^0(\C) = Lie~PGL(\alpha)$ where
$PGL(\alpha) = GL(\alpha)/\C^*(\een_{n_1},\hdots,\een_{n_k})$. 

For $\lambda
= (\lambda_1,\hdots,\lambda_k) \in \C^k$
such that $\sum_i n_i \lambda_i = 0$ we consider the element
$\underline{\lambda} = (\lambda_1 \een_{n_1},\hdots,\lambda_k \een_{n_k})$
in $M_{\alpha}^0(\C)$. The inverse image
$\mu_{\C}^{-1}(\underline{\lambda})$ is a $GL(\alpha)$-closed
affine subvariety of $rep_{\alpha}~\Q$. 

In \cite{Ginzburg:2000} V. Ginzburg proved the following coadjointness result using the
results of the preceding sections.

\begin{theorem}[Ginzburg] Assume that $\mu^{-1}_{\C}(\underline{\lambda})$
is irreducible and
that $PGL(\alpha)$ acts freely on $\mu_{\C}^{-1}(\underline{\lambda})$,
then the
quotient variety (the orbit space)
\[
\mu^{-1}_{\C}(\underline{\lambda}) / GL(\alpha) \]
is a coadjoint orbit for the necklace Lie algebra $\N_Q$.
\end{theorem}

Using results of W. Crawley-Boevey \cite{Crawley:1999} we will identify the
situations $(\alpha,\lambda)$
satisfying the conditions of the theorem. For $\lambda \in \C^k$ as above,
W. Crawley-Boevey and
M. Holland introduced and studied the {\em deformed preprojective algebra}
\[
\Pi_{\lambda} = \dfrac{\C \Q}{(m - \lambda)} \]
where $\lambda = \lambda_1 e_1 + \hdots + \lambda_k e_k \in \C \Q$. From
\cite{CrawleyHolland:1998}
we recall that $\mu^{-1}_{\C}(\underline{\lambda})$ is the scheme of
$\alpha$-dimensional representations
$\underline{rep}_{\alpha}~\Pi_{\lambda}$ of the deformed preprojective
algebra $\Pi_{\lambda}$.

We recall the characterization due to V. Kac \cite{Kac:1980} of the
dimension vectors of indecomposable
representations of the quiver $Q$.
To a vertex $v_i$ in which $Q$ has no loop, we define a {\em reflection}
\index{reflection}  $\Z^k \rTo^{r_i} \Z^k$ by
\[
r_i(\alpha) = \alpha - T_Q(\alpha,\epsilon_i) \epsilon_i \]
where $\epsilon_i = (\delta_{1i},\hdots,\delta_{ki})$.
The {\em Weyl group of the quiver $Q$} \index{Weyl group!of quiver}
\index{quiver!Weyl group}
$Weyl_Q$ is the subgroup of $GL_k(\Z)$ generated by all reflections $r_i$.

A {\em root} of the quiver $Q$ \index{root} \index{quiver!root} is a
dimension vector $\alpha \in \N^k$
such that $rep_{\alpha}~Q$ contains indecomposable representations. All
roots have connected
support. A root is said to be
\[
\begin{cases}
\text{{\em real}} &~\text{ if $\chi_Q(\alpha,\alpha) = 1$ }\\
\text{{\em imaginary}} &~\text{ if $\chi_Q(\alpha,\alpha) \leq 0$ }
\end{cases}
\]
\index{root!real} \index{root!imaginary} For a fixed quiver $Q$ we will
denote the set of all
roots, real roots and imaginary roots respectively by $\Delta, \Delta_{re}$
and $\Delta_{im}$.
With $\Pi$ we denote the set $\{ \epsilon_i~\mid~v_i~\text{has no loops}~\}$.
The {\em
fundamental set of roots} \index{roots!fundamental set} \index{fundamental
set of roots} is defined
to be the following set of dimension vectors
\[
F_Q = \{ \alpha \in \N^k - \underline{0}~\mid~T_Q(\alpha,\epsilon_i) \leq
0~\text{and $supp(\alpha)$
is connected}~\} \] Kac's result asserts
that
\[
\begin{cases}
\Delta_{re} &= Weyl_Q.\Pi \cap \N^k \\
\Delta_{im} &= Weyl_Q.F_Q \cap \N^k
\end{cases}
\]

\begin{example}
The quiver $Q$
and double quiver $\Q$ appearing in the study of Calogero phase space (see
\cite{Wilson:1998} and
\cite{Ginzburg:1999}) which stimulated the above generalizations are
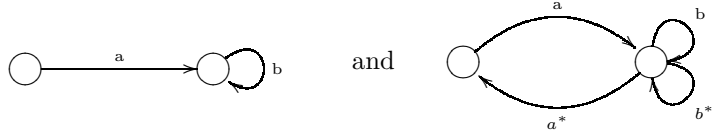
\[
\xy
\POS (0,0) *+{\txt{\tiny }}*\cir<6pt>{} ="v1" , (25,0) *+{\txt{\tiny
}}*\cir<6pt>{}="v2"
 \POS"v2" \ar@(ur,dr)^{\txt{\tiny b}}
\POS"v1" \ar^{\txt{\tiny a}} "v2"
\endxy~\qquad \text{and} \qquad
\xy
\POS (0,0) *+{\txt{\tiny }}*\cir<6pt>{} ="v1" , (25,0) *+{\txt{\tiny
}}*\cir<6pt>{}="v2"
 \POS"v2" \ar@(u,r)^{\txt{\tiny b}}
 \POS"v2" \ar@(r,d)^{\txt{\tiny $b^*$}}
\POS"v1" \ar @/^4ex/^{\txt{\tiny a}} "v2"
\POS"v2" \ar @/^4ex/^{\txt{\tiny $a^*$}} "v1"
\endxy
\]
The Euler- and Tits form of the quiver $Q$ are determined by the matrices
\[
\chi_Q = \begin{bmatrix} 1 & -1 \\ 0 & 0 \end{bmatrix} \qquad
\text{and} \qquad
T_Q = \begin{bmatrix} 2 & -1 \\ -1 & 0 \end{bmatrix}
\]
The root-system for $Q$ is easy to work out. We have
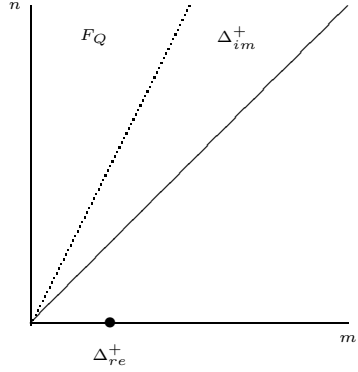
\[
\begin{xy} 0;/r.5pc/:
(0,0);(20,20) **@{-};
(0,0);(10,20) **@{.};
(5,0) *+{\bullet};
(0,0);
(20,0) **@{-};
(0,0);
(0,20) **@{-};
(4,18) *+{\txt{\tiny{$F_Q$}}};
(13,18) *+{\txt{\tiny{$\Delta^+_{im}$}}};
(20,-1) *+{\txt{\tiny $m$}};
(-1,20) *+{\txt{\tiny $n$}};
(5,-2) *+{\txt{\tiny $\Delta^+_{re}$}}
\end{xy}
\]
\[
\begin{cases}
F_Q &= \{ (m,n)~\mid~n \geq 2m \} \\
\Delta_{im}^+ &= \{ (m,n)~\mid~n \geq n \} \\
\Pi = \Delta_{re}^+ &= \{ (1,0) \}
\end{cases}
\]
\end{example}

Fix $\lambda \in \C^k$ and denote $\Delta^+_{\lambda}$ to be the set of
positive roots $\beta=(b_1,\hdots,b_k)$
for $Q$ such that $\lambda.\beta = \sum_i \lambda_i b_i = 0$. With
$S_{\lambda}$ (resp.
$\Sigma_{\lambda}$) we denote the
subsets of dimension vectors $\alpha$ which are roots for $Q$ such that
\[
1 - \chi_Q(\alpha,\alpha) \quad \geq~(\text{resp.} >) \quad r -
\chi_Q(\beta_1,\beta_1) - \hdots - \chi_Q(\beta_r,\beta_r)
\]
for all decompositions $\alpha = \beta_1 + \hdots + \beta_r$ with the
$\beta_i \in \Delta^+_{\lambda}$.
The main results of \cite{Crawley:1999} can be summarized into :

\begin{theorem}[W. Crawley-Boevey]
\begin{enumerate}
\item{$\alpha \in S_0$ if and only if $\mu_{\C}$ is a flat morphism. In
this case, $\mu_{\C}$ is also
surjective.}
\item{$\alpha \in \Sigma_{\lambda}$ if and only if $\Pi_{\lambda}$ has a
simple $\alpha$-dimensional
representation. In this case, $\mu^{-1}_{\C}(\underline{\lambda})$ is a
reduced and irreducible
complete intersection of dimension $1+\alpha.\alpha - 2
\chi_Q(\alpha,\alpha)$.}
\end{enumerate}
\end{theorem}

Using the results of \cite{LeBruynProcesi:1990} one verifies that the set of dimension
vectors of simple representations of $\Q$ coincides with the fundamental set $F_Q$. As any
simple $\Pi_{\lambda}$-representation is a simple $\Q$-representations it follows that
$\Sigma_{\lambda} \rInto F_Q$.

\begin{example}
For the Calogero-example above, we have
\begin{enumerate}
\item{The set $S_0$ consisting of all $(m,n)$ such that the complex moment
map $\mu_{\C}$ is
surjective and flat is the set of roots
\[
S_0 = \{ (m,n)~\mid~n \geq 2m-1 \} \sqcup \{ (1,0) \} \]}
\item{The set $\Sigma_0$ of dimension vectors $(m,n)$ of simple representations
of the preprojective algebra $\Pi_0$ is the set of roots
\[
\Sigma_0 = \{ (m,n)~\mid~n \geq 2m \} \sqcup \{ (1,0) \} \]
which is $F_Q \sqcup \{ (1,0) \}$.}
\item{For $\lambda = (-n,m)$ with $gcd(m,n) = 1$, the set
$\Sigma_{\lambda}$ of dimension vectors
of simple representations of the deformed preprojective algebra is the set
of roots
\[
\Sigma_{\lambda} = \{ k.(m,n)~\mid~k \in \N_+ \} \]
with unique minimal element $(m,n)$.}
\end{enumerate}
For the first two parts the essential calculation is to verify the conditions on the
decomposition $(m,n) = (m-1,n) + (1,0)$.
\end{example}

We obtain the following combinatorial description of the couples $(\alpha,\lambda)$ for
which Ginzburg's criterium applies.

\begin{theorem} $\mu^{-1}_{\C}(\underline{\lambda})$ is irreducible with a
free action of $PGL(\alpha)$
(and hence $\mu^{-1}_{\C}(\underline{\lambda})/GL(\alpha)$ is a coadjoint
orbit for $\N_Q$) if and only if
$\alpha$ is a minimal non-zero element of $\Sigma_{\lambda}$.
\end{theorem}

\begin{proof} We know that $\mu^{-1}_{\C}(\underline{\lambda}) =
\underline{rep}_{\alpha}~\Pi_{\lambda}$.
By a result of M. Artin \cite{Artin:1969} one knows that the geometric
points of the quotient scheme
$\underline{rep}_{\alpha}~\Pi_{\lambda}/GL(\alpha)$ are the isomorphism
classes of $\alpha$-dimensional
semi-simple representations of $\Pi_{\lambda}$. Moreover, the
$PGL(\alpha)$-stabilizer of a point
in $\underline{rep}_{\alpha}~\Pi_{\lambda}$ is trivial if and only if
it determines a simple
$\alpha$-dimensional representation of $\Pi_{\lambda}$. The result follows
from this and the
results recalled above.
\end{proof}

\begin{example}
Consider the special case when $\lambda = (-n,1)$ and $\alpha = (1,n)$ the
unique minimal element
in $\Sigma_{\lambda}$, then it follows from \cite{Wilson:1998} that we have
canonical identifications of the
quotient varieties
\[
iss_{\alpha}~\Pi_{\lambda} \simeq Calo_n  \]
where $Calo_n$ is the phase space of $n$ Calogero
particles. In particular, $Calo_n$ is a coadjoint
orbit.
Wilson \cite{Wilson:1998} has shown that
\[
Gr^{ad} = \bigsqcup_n~Calo_n \]
where $Gr^{ad}$ is the adelic Grassmannian which can be thought of as the space
parametrizing isomorphism classes of right ideals in the first Weyl algebra $A_1(\C) =
\C \langle x,y \rangle/(xy-yx-1)$ by \cite{CanningsHolland:1994}. In \cite{BerestWilson:1999} it
is shown that there is a non-differentiable action of the automorphism group of $A_1(\C)$ on
$Gr^{ad}$ having a transitive action on each of the $Calo_n$. It was then conjectured by
Y. Berest and G. Wilson that $Calo_n$ might be a coadjoint orbit for a central extension of the
automorphism group.
\end{example}

\begin{example}
M. Holland and W. Crawley-Boevey have a conjectural extension of the foregoing example. Let $Q'$
be an extended Dynkin quiver on $k$ vertices $\{ v_1,\hdots,v_k \}$ 
with minimal imaginary root $\delta = (d_1,\hdots,d_k)$.  A vertex $v_i$ is said to be an extending vertex
provided $d_i = 1$. Consider the quiver $Q$ on $k+1$ vertices $\{ v_0,v_1,\hdots,v_k \}$ which is
$Q'$ on the last $k$ vertices and there is one extra arrow from $v_o$ to an extending vertex $v_i$.
For a generic $\lambda' = (\lambda_1,\hdots,\lambda_k)$ they defined a noncommutative algebra
$\mathcal{O}^{\lambda'}$ extending the role of the Weyl algebra in the previous example. They conjecture
that there is a bijection between the isomorphism classes of stably free right ideals in
$\mathcal{O}^{\lambda}$ and points in
\[
\sqcup_n~\mu^{-1}_{\C}(\lambda_n)/GL(\alpha_n)
\]
where $\alpha_n = (1,n \delta)$ and $\lambda_n = (-n \lambda'.\delta,\lambda')$. This remains to be seen
but from our theorem we deduce that each of the quotient varieties
$\mu^{-1}_{\C}(\lambda_n)/GL(\alpha_n)$ is a coadjoint orbit for the necklace Lie algebra $\N_Q$.
\end{example}

If $\alpha \in \Sigma_{\lambda}$ but not minimal, there are several {\em
representation types}
$\tau = (m_1,\beta_1;\hdots,m_v,\beta_v)$ of semi-simple
$\alpha$-dimensional representations of
$\Pi_{\lambda}$ with the $\beta_i \in \Sigma_{\lambda}$ and $\sum m_i
\beta_i = \alpha$ and the
$m_i$ determine the multiplicities of the simple components. With
$iss_{\alpha}(\tau)$ we denote the
subvariety of the quotient variety $iss_{\alpha}~\Pi_{\lambda} =
rep_{\alpha}~\Pi_{\lambda}/GL(\alpha)$
consisting of all semi-simple representations of type $\tau$.

Consider the algebra $A_Q = \C[\N_Q] \otimes_{\C} \C \Q$ which has a
natural {\em trace map}
$tr : A_Q \rTo \C[\N_Q]$ mapping an oriented cycle in $\Q$ to the
corresponding necklace word and
all open paths to zero. With $Aut_Q$ we denote the automorphism group of
trace preserving
$\C$-algebra automorphisms of $A_Q$ which preserve the moment element $m =
\sum_{a \in Q_a} [a,a^*]$. A
natural extension of the above coadjoint orbit  result would be a positive
solution to the
following problem.

\begin{question} Does $Aut_Q$ act transitively on every stratum
$iss_{\alpha}(\tau)$ ?
\end{question}

\section{The smooth locus of $\Pi_0$.}

In this section and the next we try to explain why exactly the couples
$(\alpha,\lambda)$ with $\lambda.\alpha = 0$ and
$\alpha$ a minimal non-zero element of $\Sigma_{\lambda}$ give rise to
coadjoint orbits. 

The path
algebra $\C \Q$ of the double quiver $\Q$ is formally smooth in the sense of
\cite{CuntzQuillen:1995}, that is, it
has the lifting property with
respect to nilpotent ideals. Hence, $\C \Q$ is the coordinate ring of a
noncommutative affine
manifold and has
a good theory of differential forms (acyclicity). 

On the other hand, we will see that the
deformed preprojective algebras
$\Pi_{\lambda}$ are {\em never} formally smooth. 
For this reason, the
differential forms of $\C \Q$ when restricted to $\Pi_{\lambda}$ may have
rather unpredictable
behavior. 

Still, it may be possible that certain representation spaces
$rep_{\alpha}~\Pi_{\lambda}$
are smooth and we need a notion of noncommutative (formal) smoothness
depending on the
dimension vector $\alpha$. This notion is Cayley-smoothness as introduced
by C. Procesi in
\cite{Procesi:1987} and studied in detail in \cite{LeBruyn:book}.

Let $\alpha = (n_1,\hdots,n_k)$ and set $n = \sum_i n_i$. With
$\wis{alg}~{\bf @}_{\alpha}$ we denote the
category of all $V$-algebras $A$ equipped with a trace map $tr : A \rTo A$
(that is,
such that for all $a,b \in A$ we have $tr(a)b = b tr(a)$, $tr(ab) = tr(ba)$
and $tr(tr(a)b) = tr(a)tr(b)$)
satisfying $tr(1) = n$ and the formal Cayley-Hamilton identity of degree
$n$, see \cite{Procesi:1987} such that $tr(e_i) = n_i$.
Morphisms in $\wis{alg}~{\bf @}_{\alpha}$ are trace preserving algebra morphisms.
 An $\alpha$-{\em Cayley
smooth} algebra $A$ is an
algebra in $\wis{alg}~{\bf @}_{\alpha}$ having the lifting property with
respect to nilpotent ideals in
$\wis{alg}~{\bf @}_{\alpha}$. That is, every diagram
\[
\begin{diagram}
B & \rOnto^{\pi} & \frac{B}{I} \\
& \luDotsto_{\exists \tilde{\phi}} & \uTo^{\phi} \\
& & A
\end{diagram}
\]
with $B,\tfrac{B}{I}$ in $\wis{alg}~{\bf @}_{\alpha}$, $I$ a nilpotent
ideal and $\pi$ and $\phi$ trace preserving
maps, can be completed with a trace preserving algebra map $\tilde{\phi}$.
It is proved in
\cite{LeBruyn:book} that $A$ is $\alpha$-Cayley smooth if and only if the
scheme $\underline{rep}_{\alpha}~A$ of
$\alpha$-dimensional representations of $A$ is a smooth $GL(\alpha)$-variety.

In particular, if $(\lambda,\alpha)$ is such that $\lambda.\alpha = 0$ and
$\alpha$ is a minimal
non-zero vector in $\Sigma_{\lambda}$, then the level $\alpha$ approximation
$\Pi_{\lambda}~{\bf @}_{\alpha}$ (which is the ring of
$GL(\alpha)$-equivariant maps from
$\underline{rep}_{\alpha}~\Pi_{\lambda}$ to $M_n(\C)$ with the induced
trace from $M_n(\C[\underline{rep}_{\alpha}~\Pi_{\lambda})$)
is $\alpha$-Cayley smooth. In fact, $\Pi_{\lambda}~{\bf @}_{\alpha}$ is an
Azumaya algebra over the
coadjoint orbit. A neat explanation for the description of the coadjoint
orbits would be provided by a
positive solution to the following problem.

\begin{question} \label{question2} Conversely, if $\Pi_{\lambda}~{\bf
@}_{\alpha}$ is $\alpha$-Cayley smooth, does it
follow that $\alpha$ is a minimal non-zero vector in $\Sigma_{\lambda}$ ?
More generally,
does the $\alpha$-smooth locus
of
$\Pi_{\lambda}~{\bf @}_{\alpha}$, that is the locus $Sm_{\alpha}~\Pi_{\lambda}$ in $\underline{iss}_{\alpha}~\Pi_{\lambda}$
such that $\underline{rep}_{\alpha}~\Pi_{\alpha}$ is smooth along $\pi^{-1}(Sm~\Pi_{\lambda})$,
coincide with the Azumaya locus ?
\end{question}

We will give an affirmative solution in the special case of the preprojective
algebra $\Pi_0$. By a result of 
W. Crawley-Boevey \cite{Crawley:1999b}, we can control the $Ext^1$-spaces of representations of $\Pi_0$.
Let $V$ and $W$ be representations of $\Pi_0$ of
dimension vectors $\alpha$ and $\beta$, then we have
\[
dim_{\C}~Ext^1_{\Pi_0}(V,W) = dim_{\C}~Hom_{\Pi_0}(V,W) + dim_{\C}~Hom_{\Pi_0}(W,V) - T_Q(\alpha,\beta)
\]
For $\xi \in \underline{iss}_{\alpha}~\Pi_0$ to belong to the
smooth locus $\xi \in Sm_{\alpha}~\Pi_0$ it is necessary and
sufficient that $\underline{rep}_{\alpha}~\Pi_0$ is smooth along the orbit $\Oscr(M_{\xi})$ where $M_{\xi}$
is the semi-simple $\alpha$-dimensional representation of $\Pi_0$ corresponding to $\xi$. 

Assume that
$\xi$ is of type $\tau = (e_1,\alpha_1;\hdots;e_z,\alpha_z)$, that is,
\[
M_{\xi} = S_1^{\oplus e_1} \oplus \hdots \oplus S_z^{\oplus e_z} \]
with $S_i$ a simple $\Pi_0$-representation of dimension vector $\alpha_i$. Then, the normal space
to the orbit $\Oscr(M_{\xi})$ is determined by $Ext^1_{\Pi_o}(M_{\xi},M_{\xi})$ and can be
depicted by a local quiver setting $(Q_{\xi},\alpha_{\xi})$ where $Q_{\xi}$ is a quiver
on $z$ vertices having as many arrows from vertex $i$ to vertex $j$ as the dimension
of $Ext^1_{\Pi_0}(S_i,S_j)$ and where $\alpha_{\xi} = \alpha_{\tau} = (e_1,\hdots,e_z)$. 
Applying the Luna slice theorem \cite{Luna:1973} we have

\begin{lemma} With notations as above, $\xi \in Sm_{\alpha}~\Pi_0$ if and only if
\[
dim~GL(\alpha) \times^{GL(\alpha_{\xi})} Ext^1_{\Pi_0}(M_{\xi},M_{\xi}) =
 dim_{M_{\xi}}~\underline{rep}_{\alpha}~\Pi_0 \]
\end{lemma}

As we have enough information to compute both sides, we can prove :

\begin{theorem} If $\xi \in \underline{iss}_{\alpha}~\Pi_0$ with
$\alpha = (a_1,\hdots,a_k) \in S_{\underline{0}}$, then $\xi \in Sm_{\alpha}~\Pi_0$
if and only if $M_{\xi}$ is a simple $n$-dimensional representation of $\Pi_0$. That is, the smooth locus of $\Pi_0$
coincides with the Azumaya locus.
\end{theorem}

\begin{proof}
Assume that $\xi$ is a point of semi-simple representation type $\tau = (e_1,\alpha_1;\hdots;e_z,\alpha_z)$, that is,
\[
M_{\xi} = S_1^{\oplus e_1} \oplus \hdots \oplus S_z^{\oplus e_z} \qquad \text{with} \qquad
dim(S_i) = \alpha_i \]
and $S_i$ a simple $\Pi_0$-representation. 
We have
\[
\begin{cases}
dim_{\C}~Ext^1_{\Pi_0}(S_i,S_j) &= -T_Q(\alpha_i,\alpha_j) \qquad i \not= j \\
dim_{\C}~Ext^1_{\Pi_0}(S_i,S_i) &= 2 - T_Q(\alpha_i,\alpha_i)
\end{cases}
\]
But then, the dimension of $Ext^1_{\Pi_0}(M_{\xi},M_{\xi})$ is equal to
\[
\sum_{i=1}^z (2- T_Q(\alpha_i,\alpha_i))e_i^2 + \sum_{i \not= j} e_ie_j(-T_Q(\alpha_i,\alpha_j) =
2 \sum_{i=1}^z e_i - T_Q(\alpha,\alpha) \]
from which it follows immediately that
\[
dim~GL(\alpha) \times^{GL(\alpha_{\xi})} Ext^1_{\Pi_0}(M_{\xi},M_{\xi}) = \alpha.\alpha +
\sum_{i=1}^z e_i^2 - T_Q(\alpha,\alpha) \]
On the other hand, as $\alpha \in S_{\underline{0}}$ we know
that
\[
dim~\underline{rep}_{\alpha}~\Pi_0 = \alpha.\alpha - 1 + 2 p_Q(\alpha) = 
\alpha.\alpha - 1 + 2 - 2 \chi_Q(\alpha,\alpha) = \alpha.\alpha + 1 - T_Q(\alpha,\alpha) \]
But then, equality occurs if and only if $\sum_i e_i^2 = 1$, that is, $\tau = (1,\alpha)$ or
$M_{\xi}$ is a simple $n$-dimensional representation of $\Pi_0$.
\end{proof}

In particular it follows that the preprojective algebra $\Pi_0$ is {\em never} formally smooth
as this implies that all the representation varieties must be smooth.
Further, as $\vec{v_i} = (0,\hdots,1,0,\hdots,0)$ are dimension vectors of simple
representations of $\Pi_0$ it follows that $\Pi_0$ is $\alpha$-smooth if and only if
$\alpha = \vec{v_i}$ for some $i$. 

\begin{example}  Let $Q$ be an extended Dynkin diagram and $\delta$ the minimal imaginary
root, then $\delta \in S_{0}$. The dimension of the quotient variety
\[
\begin{split}
dim~iss_{\delta}~\Pi_0 &= dim~rep_{\delta}~\Pi_0 - \delta.\delta + 1 \\
&= 2
\end{split}
\]
so it is a surface. The only other semi-simple $\delta$-dimensional representation of $\Pi_0$ is
the trivial representation. By the theorem, this must be an
isolated singular point of $iss_{\delta}~Q$. In fact, one can show that $iss_{\delta}~\Pi_0$ is
the Kleinian singularity corresponding to the extended Dynkin diagram $Q$.
\end{example}

\section{A sheaf of algebras.}

We will prove that $\alpha$-Cayley smoothness of a closely related sheaf of
algebras is equivalent
to $\alpha$ being a minimal non-zero vector of $\Sigma_{\lambda}$. Recall
that $rep_{\alpha}~\Q$
admits  a hyper-K\"ahler structure (that is, an action of the quaternion
algebra
$\mathbb{H} = \R.1 \oplus \R.i \oplus \R.j \oplus \R.k$) defined for all
arrows $a \in Q_a$ and all
arrows $b \in \Q_a$ by the formulae, see for example \cite{Crawley:1999b}
\[
\begin{split}
(i.V)_b &= iV_b \\
(j.V)_a &= -V_{a^*}^{\dagger} \quad (j.V)_{a^*} = V_a^{\dagger} \\
(k.V)_a &= -iV_{a^*}^{\dagger} \quad (k.V)_{a^*} = i V_a^{\dagger}
\end{split}
\]
where this time we denote the Hermitian adjoint of a matrix $M$ by
$M^{\dagger}$ to distinguish
it from the star-operation on the arrows of the double quiver $\Q$. Let
$U(\alpha)$ be the product
of unitary groups $U_{n_1} \times \hdots \times U_{n_k}$ and consider the
{\em real moment map}
\[
rep_{\alpha}~\Q \rTo^{\mu_{\R}} Lie~U(\alpha)~\qquad~V \mapsto
\underset{\underset{b \in \Q_a}{\xy /r.1pc/:
\POS (15,0) *\cir<4pt>{}*+{\txt\tiny{}} ="v2"
   , (0,0) *\cir<4pt>{}*+{\txt\tiny{}} ="v3"
\POS"v2" \ar "v3"_{b} \endxy}}{\sum}~\frac{i}{2}[V_b,V_b^{\dagger}] \]
For $\lambda \in \R^k$, multiplication by the quaternion-element $h =
\frac{i+k}{\sqrt{2}}$ gives
a homeomorphism between the real varieties
\[
\mu^{-1}_{\C}(\underline{\lambda}) \cap \mu^{-1}_{\R}(\underline{0}) \rTo^{h.}
\mu_{\C}^{-1}(\underline{0}) \cap
\mu^{-1}_{\R}(i \underline{\lambda}) \]
Moreover, the hyper-K\"ahler structure commutes with the base-change action
of $U(\alpha)$, whence we have
a natural one-to-one correspondence between the quotient spaces
\[
(\mu^{-1}_{\C}(\underline{\lambda}) \cap
\mu^{-1}_{\R}(\underline{0}))/U(\alpha) \rTo^{h.}
(\mu_{\C}^{-1}(\underline{0}) \cap
\mu^{-1}_{\R}(i \underline{\lambda}))/U(\alpha) \]
see \cite{Crawley:1999b} for more details. By results of Kempf and Ness
\cite{KempfNess:1979} we can
identify the left hand side as the quotient variety
$iss_{\alpha}~\Pi_{\lambda}$ and by results of
A. King \cite{King:1994} we can identify the right hand side as the moduli
space
$M_{\alpha}^{ss}(\Pi_0,\lambda)$ of $\lambda$-semistable
$\alpha$-dimensional representations of the
preprojective algebra $\Pi_0$. Recall that a representation $V \in
rep_{\alpha}~\Q$ is said to be
$\lambda$-(semi)stable if and only if for every proper subrepresentation
$W$ of $V$ say with
dimension vector $\beta$ we have $\lambda.\beta > 0$ (resp. $\lambda.\beta
\geq 0$). The scheme
$\underline{rep}_{\alpha}^{ss}(\Pi_0,\lambda)$ of $\lambda$-semistable
$\alpha$-dimensional representations
of $\Pi_0$ is the intersection of $\mu_{\C}^{-1}(\underline{0})$ with the
subvariety of $\lambda$-semistable
representations in $rep_{\alpha}~\Q$. The corresponding moduli space
$M_{\alpha}^{ss}(\Pi_0,\lambda)$
classifies isomorphism classes of direct sums of $\lambda$-stable
representations of
$\Pi_0$ of total dimension $\alpha$. In view of the explicit form of the
hyper-K\"ahler structure it
follows that the deformed
preprojective algebra $\Pi_{\lambda}$ has semi-simple representations of
dimension vector $\alpha$ of
representation type $\tau = (e_1,\beta_1;\hdots;e_r,\beta_r)$ if and only
if the preprojective
algebra $\Pi_0$ has $\lambda$-stable representations of dimension vectors
$\beta_i$ for all $1 \leq i \leq r$.
In particular, $\Pi_{\theta}$ has a simple representation of dimension
vector $\alpha$ if and only if
$\Pi_0$ has a $\theta$-stable representation of dimension vector $\alpha$.

Taking locally the algebras of $GL(\alpha)$-equivariant maps from
$\underline{rep}_{\alpha}^{ss}(\Pi_0,\lambda)$
to $M_n(\C)$ defines a sheaf of algebras in $\wis{alg}~{\bf @}_{\alpha}$,
$\mathcal{A}_{\lambda,\alpha}$ on the moduli space
$M_{\alpha}^{ss}(\Pi_0,\lambda)$. The main result of this section asserts
the following.

\begin{theorem} With notations as above, for $\alpha \in \Sigma_{\lambda}$
the following are equivalent :
\begin{enumerate}
\item{$\mathcal{A}_{\lambda,\alpha}$ is a sheaf of $\alpha$-Cayley smooth
algebras on the
moduli space $M_{\alpha}^{ss}(\Pi_0,\lambda)$.}
\item{$\alpha$ is a minimal non-zero vector in $\Sigma_{\lambda}$ (and
hence the quotient variety
$iss_{\alpha}~\Pi_{\lambda}$ is a coadjoint orbit for the necklace Lie
algebra $\N_Q$).}
\end{enumerate}
\end{theorem}

\begin{proof}
As $\alpha \in \Sigma_{\lambda}$ we know that $iss_{\alpha}~\Pi_{\lambda}$
has dimension $
1+\alpha.\alpha-2 \chi_Q(\alpha,\alpha) - dim~PGL(\alpha)$ which is equal to
$2 - T_Q(\alpha,\alpha)$. By the hyper-K\"ahler correspondence so is the
dimension of
$M_{\alpha}^{ss}(\Pi_0,\lambda)$, whence the open subset of
$\mu^{-1}_{\C}(\underline{0})$
consisting of $\lambda$-semistable representations has dimension
\[
1 + \alpha.\alpha - 2 \chi_Q(\alpha,\alpha) \]
as there are $\lambda$-stable representations in it (again via the
hyper-K\"ahler
correspondence). Take a $GL(\alpha)$-closed orbit $\Oscr(V)$ in this
open set. That is, $V$ is the direct sum of $\lambda$-stable subrepresentations
\[
V = S_1^{\oplus e_1} \oplus \hdots \oplus S_r^{\oplus e_r} \]
with $S_i$ a $\lambda$-stable representation of $\Pi_0$ of dimension vector
$\beta_i$ occurring in $V$
with multiplicity $e_i$ whence $\alpha = \sum_i e_i \beta_i$.

Again, the normal space in $V$ to $\Oscr(V)$ can be identified with
$Ext^1_{\Pi_0}(V,V)$.
As all $S_i$ are $\Pi_0$-representations we can determine this space by the
knowledge of
all $Ext^1_{\Pi_0}(S_i,S_j)$. 
\[
Ext^1_{\Pi_0}(S_i,S_j) = 2 \delta_{ij} - T_Q(\beta_i,\beta_j) \]
But then the dimension of the normal space to the orbit is
\[
dim~Ext^1_{\Pi_0}(V,V) = 2 \sum_{i=1}^r e_i - T_Q(\alpha,\alpha) \]
By the Luna slice theorem \cite{Luna:1973}, the \'etale local structure in
the smooth point $V$
is of the form $GL(\alpha) \times^{GL(\tau)}
Ext^1(V,V)$ where $\tau = (e_1,\hdots,e_r)$ and is therefore of dimension
\[
\alpha. \alpha + \sum_{i=1}^2 e_i^2 - T_Q(\alpha,\alpha) \]
This number must be equal to the dimension of the subvariety of
$\lambda$-semistable representations
of $\Pi_0$ which has dimension $1 + \alpha.\alpha - T_Q(\alpha,\alpha)$ if
and only if
$r=1$ and $e_1 = 1$, that is if and only if $V$ is $\lambda$-stable.
Hence, if $rep_{\alpha}^{ss}(\Pi_0,\lambda)$ is smooth, then $\alpha$ must
be a minimal non-zero
vector in the set of dimension vectors of $\lambda$-stable representations
of $\Pi_0$ and hence by
the hyper-K\"ahler correspondence, $\alpha$ is a minimal non-zero vector in
$\Sigma_{\lambda}$.

Conversely, if $\alpha$ is a minimal vector in $\Sigma_{\lambda}$, then
$iss_{\alpha}~\Pi_{\lambda}$
is a coadjoint orbit, whence smooth and hence so is
$M_{\alpha}^{ss}(\Pi_0,\lambda)$ by the
correspondence. Moreover, all $\alpha$-dimensional $\lambda$-semistable
representations must be
$\lambda$-stable by the minimality assumption and so
$rep_{\alpha}^{ss}(\Pi_0,\lambda)$ is
a principal $PGL(\alpha)$-fibration over $M_{\alpha}^{ss}(\Pi_0,\lambda)$
whence smooth. Therefore,
$\mathcal{A}_{\lambda,\alpha}$ is a sheaf of $\alpha$-Cayley smooth algebras.
\end{proof}

Question~\ref{question2} can be proved as in the case of $\Pi_0$ provided
we know that
\[
\begin{split}
dim_{\C}~Ext^1_{\Pi_{\lambda}}(M,N) &= dim_{\C}~Hom_{\Pi_{\lambda}}(M,N) +
dim_{\C}~Hom_{\Pi_{\lambda}}(N,M) \\
&- T_Q(dim~M,dim~N)
\end{split}
\]
for all (simple) $\Pi_{\lambda}$ representations $M$ and $N$. Still, we can prove that
$\Pi_{\lambda}$ can never be formally smooth and even that certain representation varieties
are not smooth.

\begin{proposition} Let $\alpha \in \Sigma_{\lambda}$ such that $2 \alpha \in \Sigma_{\lambda}$.
Then, $\underline{rep}_{2 \alpha}~\Pi_{\lambda}$ is not smooth. In particular, $\Pi_{\lambda}$ is
not formally smooth.
\end{proposition}

\begin{proof} As $\alpha \in \Sigma_{\lambda}$ we know that the local quiver $Q_{\xi}$ in
a simple representation $S$ corresponding to
 $\xi$ is a one vertex quiver having $2-T_Q(\alpha,\alpha)$ loops. That is,
\[ dim~Ext^1_{\Pi_{\lambda}}(S,S) = 2 - T_Q(\alpha,\alpha) \]
 But then, for $\xi \in \underline{iss}_{2 \alpha}~\Pi_{\lambda}$ a point corresponding to
 $S \oplus S$, the local quiver is still $Q_{\xi}$ but this time the local dimension vector
 $\alpha_{\xi} = 2$. If $\xi$ lies in the smooth locus, then by the Luna slice theorem we must have
 \[
 dim~GL(2 \alpha) \times^{GL_2} rep_{\alpha_{\xi}}~Q_{\xi} = dim~rep_{2 \alpha}~\Pi_{\lambda} \]
 The left hand side is $4 \alpha.\alpha +4 - 4 T_Q(\alpha,\alpha)$ whereas the right hand side is
 equal to (because $2 \alpha \in \Sigma_{\lambda}$) $4 \alpha.\alpha + 1 - 4 T_Q(\alpha,\alpha)$,
 a contradiction.
 \end{proof}

\end{document}